\documentclass[a4paper, 10pt, reqno]{amsart}
\usepackage[margin=1in]{geometry}
\usepackage{amsmath}
\usepackage{amssymb}
\usepackage{mathtools}
\usepackage{tikz}
\usepackage{enumerate}
\usepackage{gensymb}
\usepackage{graphicx}
\usepackage{pgfplots,array}
\usepackage{algorithm}
\usepackage{algpseudocode}
\usepackage{booktabs}
\usepackage{subcaption}
\usepackage[width=\textwidth]{caption}
\makeatletter
\g@addto@macro\bfseries{\boldmath}
\makeatother
\algnewcommand\algorithmicforeach{\textbf{for each}}
\algdef{S}[FOR]{ForEach}[1]{\algorithmicforeach\ #1\ \algorithmicdo}
\title{Optimal Gap Sequences in Shellsort for $n\leq16$ Elements}
\author{Ying Wai Lee}
\date{\today}
\begin{document}
\maketitle
\begin{abstract}
Optimal gap sequences in Shellsort, defined as gap sequences having the minimised maximum number of comparisons for a fixed number of pairwise distinct elements, are found by minimax search in reduced permutational spaces, namely Bad $(s,1)$-sorted permutations. Exact optimal gap sequences in Shellsort for $n\leq16$ pairwise distinct elements are established, and the best known gap sequences for $17\leq n\leq 30$ are listed with conjectures made. It notably discovers some optimal gap sequences consist of increments larger than the half of the total number of the elements to sort.
\end{abstract}

\section{Introduction}
Shellsort, invented by Shell (1959), is one of the most traditional sorting algorithms, in which a permutation is sorted by several iterations of linear insertion sort on the elements that stay apart by a positive integer from a gap sequence, called an increment that progressively decreases to 1. The universally optimal gap sequence to sort any number of elements by Shellsort, and the optimal asymptotic growth on the average number of comparison required in by Shellsort still remain as open problems. The currently best known gap sequence in Shellsort for large number of elements, in which the $k$-th increment $h_k$, is given by (Lee, 2021),
\begin{align}\label{Lee}
h_k=\left\lceil
\frac{\gamma^k-1}{\gamma-1}
\right\rceil
\end{align}
, where $\gamma=2.243609061420001...$, $k\in\mathbb{N}_1$ and $\lceil x\rceil$ is ceiling function. This gap sequences is known as a $\gamma$-seqence, and the first few increments of the $\gamma$-sequence given by (\ref{Lee}) are
\begin{align}
1,\, 4,\, 9,\, 20,\, 45,\, 102,\, 230,\, 516,\,1158,\,2599,\,5831,\,13082,\,29351,\,65853,\, 147748,\,331490,\,743735,\, ...\end{align}

For any $n\in\mathbb{N}_1$, one of the possible methods to find the optimal gap sequence to sort $n$ pairwise distinct elements by Shellsort, is to generate all permutations of $n$ pairwise distinct elements, all gap sequences to sort those element, and find the gap sequence that has the least maximum number of comparisons required in sorting all those permutations. It is achieved, by the definition of the optimal gap sequence, by a minimax search algorithm and further simplified in reduced permutational space.

\section{Main Result}
Tables \ref{tab:my-table_1} and 2 respectively summarise the optimal gap sequences and best known gap sequences found to sort up to 16 and 30 pairwise distinct elements respectively. For any $n=1,2,3,4,5$, the optimal gap sequence consists of 1 as the only increment, that makes optimised Shellsort in sorting elements up to 5 pairwise distinct elements being the same as the linear insertion sort. For any $n=6,7,8,10,12,16$, the optimal gap sequence consists of at least increment which is larger than $n/2$.
\clearpage\vspace*{\fill}\begin{table}[H]
\centering
\begin{tabular}{>{\centering\arraybackslash}p{.75cm}|>{\centering\arraybackslash}p{6cm}|>{\centering\arraybackslash}p{2.5cm}|>{\centering\arraybackslash}p{5cm}}
$n=$ & Optimal Gap Sequence ($s^{**}_n=$) & Index $i(s^{**}_n)=$ & Minimax Comparison ($c_n=$) \\ \hline
 1 &             &    & 0 \\
 2 & 1           &  1 & 1 \\
 3 & 1           &  1 & 3 \\
 4 & 1           &  1 & 6 \\
 5 & 1           &  1 & 10 \\
 6 & 1, 4        &  5 & 14 \\
 7 & 1, 4, 6     & 21 & 18 \\
 8 & 1, 5, 7     & 41 & 23 \\
 9 & 1, 3, 4     &  7 & 29 \\
10 & 1, 6, 9     &145 & 35 \\
11 & 1, 4, 5     & 13 & 41 \\
12 & 1, 3, 7, 11 &543& 48 \\
13 & 1, 3, 4     &  7 & 56 \\
14 & 1, 3, 4     &  7 & 64 \\
15 & 1, 3, 7     & 35 & 71 \\
16 & 1, 4, 7, 9  &165 & 78
\end{tabular}
\caption{Optimal Gap Sequences in Shellsort for $n$ Pairwise Distinct Elements}
\label{tab:my-table_1}
\end{table}
\vspace*{\fill}
\begin{table}[H]
\centering
\begin{tabular}{>{\centering\arraybackslash}p{.75cm}|>{\centering\arraybackslash}p{6cm}|>{\centering\arraybackslash}p{2.5cm}|>{\centering\arraybackslash}p{5cm}}
$n=$ & Best Known Gap Sequence ($s'_n=$) & Index $i(s'_n)=$ & Upper Bound of $c_n$ ($c_n\leq$) \\ \hline
17 & 1, 3, 4 &   7 & $87$ \\
18 & 1, 2, 5 &  10 & $98$ \\
19 & 1, 3, 5 &  11 & $105$ \\
20 & 1, 3, 4 &   7 & $117$ \\
21 & 1, 3, 4 &   7 & $126$ \\
22 & 1, 2, 3 &   4 & $157$ \\
23 & 1, 2, 3 &   4 & $173$ \\
24 & 1, 2, 3 &   4 & $183$ \\
25 & 1, 2, 3 &   4 & $195$ \\
26 & 1, 2, 3 &   4 & $219$ \\
27 & 1, 2, 3 &   4 & $230$ \\
28 & 1, 2, 3 &   4 & $243$ \\
29 & 1, 2, 3 &   4 & $263$ \\
30 & 1, 2, 3 &   4 & $275$ 
\end{tabular}
\caption{Best Known Gap Sequences in Shellsort for $n$  Pairwise Distinct Elements}
\label{tab:my-table_1b}
\end{table}\vspace*{\fill}

\clearpage\section{Methodology \& Terminology}
Let $n\in\mathbb{N}_1$ be a positive integer and $P_n$ be the set of all permutations of $\{1,2,...,n\}$. A finite set of positive integers $s\subset\mathbb{N}_1$ is said to be a gap sequence for $n$, if $1\in s$ and $\max{s}\leq n-1$. The set of all gap sequences for $n$ is deonted by $S_n$. For any $s\in S_n$ and $p\in P_n$, the number to comparisons required to sort $p$ by $s$ in Shellsort is denoted by $n_{p,s}$. Hence the minimised maximum number of comparisons, denoted by $c_n$, to sort $n$ pairwise distinct elements by Shellsort is given by \begin{align}
c_n=\min{\{\max{\left\{n_{p,s}:p\in P_n\right\}}:s\in S_n\}}
\end{align}

$s^*\in S_n$ is said to be an optimal gap sequence for $n$ if it has the best worse-case complexity in making comparison, or equivalently $\max{\{n_{p,s^*}:p\in P_n\}}=c_n$. Gap sequences can be parameterised by an index $i\in\mathbb{N}_1$, for which a gap sequence $s_i$ is given by writing $i$ in binary system and
\begin{align}
i&=\sum_{k=0}^{\lfloor\log_2{i}\rfloor}2^kb(i,k)\\
s_i&=\{1\}\cup\{k\in\mathbb{N}_2:b(i-1,k-1)=1\}
\end{align}

This method provides a systemic way to generate all gap sequences for a minimax search algorithm to find the number $c_n$. Table \ref{tab:my-table_3} summaries the first few examples of $s_i$, notice that the definition of $s_i$ is independent of the value of $n$. Every gap sequence $s\in S_n$ has a unique index $i(s)\in\mathbb{N}_1$ such that $s=s_{i(s)}$, which can be computed and given by
\begin{align}
i(s)=\frac{1}{2}+\sum_{k\in s}2^{k-2}
\end{align}

$s^{**}_n\in S_n$ is said to be the optimal gap sequence for $n$ if for any $s^*$, an optimal gap sequence for $n$, $i({s^{**}_n})\leq i({s^*})$. A minimax search algorithm is hence obtained to find the optimal gap sequence by the parameterisation index $i$, starting from $i=1$ up to $i=2^{n-2}$. By the definition of $s^{**}_n$, the optimal gap sequence must have the smallest largest increment possible among all optimal gap sequences.

The number of permutations for $n$ pairwise distinct elements is given by $n!$, $n$ factorial. Thus $|P_n|$, the number of elements in $P_n$, is $n!$. However, for any $s\in S_n$, it is not necessary to consider every element $p$ in $P_n$ to find the value of $\max{\{n_{p,s}:p\in P_n\}}$, as there are many redundancies to omit and the maximum must happen in a proper subset of $P_n$.

\vfill\begin{table}[h!]
\centering
\begin{tabular}{>{\centering\arraybackslash}p{.75cm}|p{2.5cm}||>{\centering\arraybackslash}p{.75cm}|p{2.5cm}||>{\centering\arraybackslash}p{.75cm}|p{2.5cm}}
$i=$ & $s_i=$ & $i=$ & $s_i=$ & $i=$ & $s_i=$ \\ \hline
 1 & 1        &  7 & 1, 3, 4     & 13 & 1, 4, 5 \\
 2 & 1, 2     &  8 & 1, 2, 3, 4  & 14 & 1, 2, 4, 5 \\
 3 & 1, 3     &  9 & 1, 5        & 15 & 1, 3, 4, 5 \\
 4 & 1, 2, 3  & 10 & 1, 2, 5     & 16 & 1, 2, 3, 4, 5 \\
 5 & 1, 4     & 11 & 1, 3, 5     & 17 & 1, 6 \\
 6 & 1, 2, 4  & 12 & 1, 2, 3, 5  & 18 & 1, 2, 6       
\end{tabular}\caption{First Few Examples of $s_i$}
\label{tab:my-table_3}
\end{table}

\clearpage\section{Reduced Permutational Space}
For example, $s=\{1,4\}\in S_n$ is a gap sequence of $n=16$. Consider the permutations $p_1,p_2\in P_n$.

\vfill\begin{figure}[h]
\begin{subfigure}[b]{0.49\linewidth}
\begin{tikzpicture}
\centering
\begin{axis} [width=1.2\textwidth,title=\text{$p_1=(16,12,8,4,15,11,7,3,14,10,6,2,13,9,5,1)$}, ybar,bar width=11pt,clip=false,axis y line=none, xticklabels=\empty,axis x line*=bottom,nodes near coords,ymin=0,xmajorticks=false,]
\addplot [fill=none,ytick=data]coordinates { 
 ( 1,16) 
 ( 2,12) 
 ( 3, 8) 
 ( 4, 4)
 ( 5,15) 
 ( 6,11) 
 ( 7, 7) 
 ( 8, 3)
 ( 9,14) 
 (10,10) 
 (11, 6) 
 (12, 2)
 (13,13) 
 (14, 9) 
 (15, 5) 
 (16, 1)
};
\end{axis} 
\end{tikzpicture}
\end{subfigure}\hfill
\begin{subfigure}[b]{0.49\linewidth}
\begin{tikzpicture}
\centering
\begin{axis} [width=1.2\textwidth,title=\text{$p_2=(13,9,5,1,14,10,6,2,15,11,7,3,16,12,8,4)$}, ybar,bar width=11pt,clip=false,axis y line=none, xticklabels=\empty,axis x line*=bottom,nodes near coords,ymin=0,xmajorticks=false,]
\addplot [fill=none,ytick=data]coordinates { 
 ( 1,13) 
 ( 2, 9) 
 ( 3, 5) 
 ( 4, 1)
 ( 5,14) 
 ( 6,10) 
 ( 7, 6) 
 ( 8, 2)
 ( 9,15) 
 (10,11) 
 (11, 7) 
 (12, 3)
 (13,16) 
 (14,12) 
 (15, 8) 
 (16, 4)
};
\end{axis} 
\end{tikzpicture}\end{subfigure}\caption{Illustrations of Permutations {$p_1\in P_{n,(s,1)}$ and $p_2\in P_n\setminus P_{n,(s,1)}$}}\end{figure}
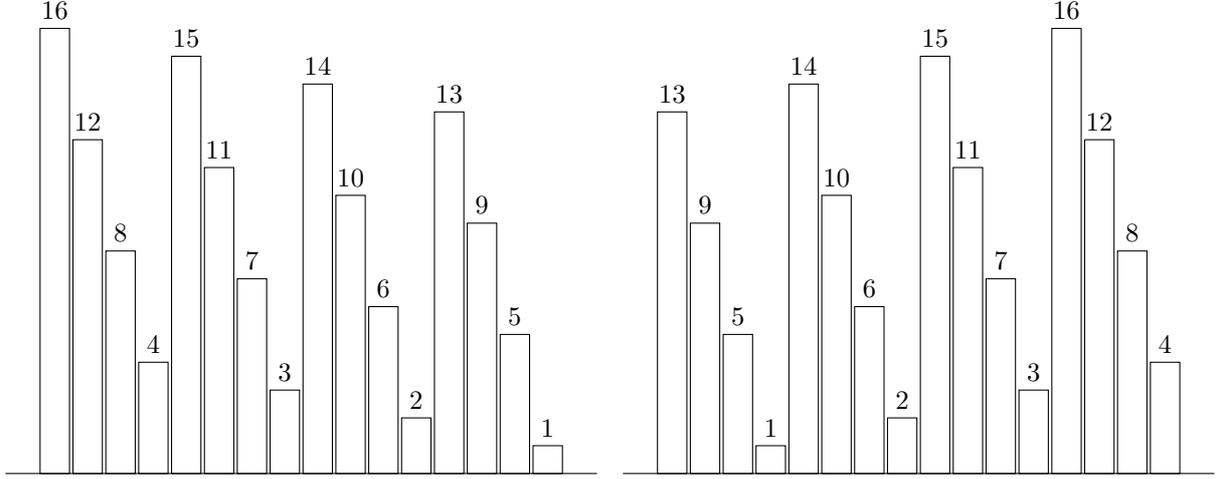

After the first pass, which is 4 gapped linear insertion sorts in ascending order for $4\in s$, $p_1$ becomes the 4-sorted $p_2$, and $p_2$ remains unchanged. The numbers of comparisons of $p_1$ and $p_2$ for the increment $4\in s$ are $(1+2+3)\times 4=24$ and $3\times 4=12$ respectively. Hence, it suffices to consider a proper subset of $P_n$ for those permutations $p$ that behave similar to $p_1$ under the gap sequence $s$.

Let $n\in\mathbb{N}_1$, $s\in S_n$ and $p\in P_n$. $p$ is said to be Bad $(s,1)$-sorted, if every sub-sequence of every $s(1)$, the largest increment in $s$, elements apart being the worse case in linear insertion sort in the sense of maximum number of comparisons. The set of all Bad $(s,1)$-sorted permutations is denoted by $P_{n,(s,1)}$. Thus, $c_n$ can be rewritten as the minimum of the maximum on a proper subset of domain, and \begin{align}
c_n=\min{\left\{\max{\left\{n_{p,s}:p\in P_{n,(s,1)}\right\}}:s\in S_n\right\}}
\end{align}

For any $n\in\mathbb{N}_1$ and $s\in S_n$, $P_{n,(s,1)}$ is always a proper subset of $P_n$. In the example, $n=16$, $s=\{4,1\}\in S_n$, $p_1\in P_{n,(s,1)}$ and $p_2\in P_n\setminus P_{n,(s,1)}$. The number of permutations that sufficient to consider for $s$ is largely reduced from $|P_n|=16!= 20\,922\,789\,888\,000$ to $|P_{n,(s,1)}|=16!/(4!)^4=63\,063\,000$. Table \ref{tab:my-table_2} summaries all of the cases of $s\in S_n$ for $n=16$.

\vfill\begin{table}[h!]
\centering
\begin{tabular}{>{\centering\arraybackslash}p{1.1cm}|>{\centering\arraybackslash}p{3.3cm}||>{\centering\arraybackslash}p{1.1cm}|>{\centering\arraybackslash}p{3.3cm}||>{\centering\arraybackslash}p{1.1cm}|>{\centering\arraybackslash}p{3.3cm}}
$s(1)=$ & $\left|P_{n,(s,1)}\right|=$ & $s(1)=$ & $\left|P_{n,(s,1)}\right|=$ & $s(1)=$ & $\left|P_{n,(s,1)}\right|=$ \\ \hline
 1 &           1 &  6 &   4 036 032 000 & 11&    653 837 184 000 \\
 2 &      12 870 &  7 &  18 162 144 000 & 12&  1 307 674 368 000 \\
 3 &   2 018 016 &  8 &  81 729 648 000 & 13&  2 615 348 736 000 \\
 4 &  63 063 000 &  9 & 163 459 296 000 & 14&  5 230 697 472 000 \\
 5 & 672 672 000 & 10 & 326 918 592 000 & 15& 10 461 394 944 000
\end{tabular}\caption{Numbers of Elements in $P_{n,(s,1)}$ for $n=16$}\label{tab:my-table_2}\end{table}

\clearpage\section{Implementation of Minimax Search}
Adopted from the pseudocode, minimax searchs are performed. The full execution results for $n\leq 16$, in Tables \ref{table:12} and \ref{table:16}, and the partial execution results for $17\leq n\leq 30$, in Tables \ref{table:20} and \ref{table:30}, are listed. 
\hfill \\
\vfill
\setlength{\intextsep}{0pt} \begin{algorithm}
  \caption{Minimax Algorithm in Finding the Optimal Gap Sequence in Shellsort}
  \begin{algorithmic}
    \Procedure{Minimax}{$n$}
    \State $c_n \gets n(n-1)/2$ 
    \ForEach {$i \in \{1,2,...,2^{n-2}\} $}
    \State $s \gets s_i$ 
    \State $n_i \gets 0$
    \ForEach {$p \in P_{n,(s,1)} $}
    \State $t_{s,p}\gets$ number of comparisons to sort $p$ by $s$ in Shellsort
    \State $n_i \gets \max{(n_i,t_{s,p})}$
    \If{$n_i\geq c_n$}
    \State \textbf{break}
    \EndIf
    \EndFor
    \If{$n_i< c_n$}
    \State $c_n \gets n_i$
    \State \Return $i$, $s_i$, $n_i$
    \EndIf
    \EndFor
    \EndProcedure
  \end{algorithmic}
\end{algorithm}

\vfill For $n=1,2,3,4,5$, the minimax search yields the optimal gap sequence to be exactly the same as the linear insertion sort, and Shellsort starts to require lesser numbers of comparisons for $n\geq6$. Although the implementations of minimax search algorithm and reduced permutational space save a lot of computation power, the fast growth of the number of permutations required to consider still makes the search to last for a long time of computation to terminate. Due to the limitation of efficient computational power and the lack of further techniques in simplification, the optimal gap sequences for $n\geq17$ are not yet established. 

Table \ref{tab:Shell_Linear} summaries the known improvements for the numbers of comparisons of Shellsort by the optimal gap sequences over the linear insertion sort in respective worst cases, and the improvement is more pronounced as $n\in\mathbb{N}_1$ increases.
\vfill\begin{table}[H]
\begin{tabular}{>{\centering\arraybackslash}p{2.25cm}|
>{\centering\arraybackslash}p{0.475cm}>{\centering\arraybackslash}p{0.475cm}>{\centering\arraybackslash}p{0.475cm}>{\centering\arraybackslash}p{0.475cm}>{\centering\arraybackslash}p{0.475cm}>{\centering\arraybackslash}p{0.475cm}>{\centering\arraybackslash}p{0.475cm}>{\centering\arraybackslash}p{0.475cm}>{\centering\arraybackslash}p{0.475cm}>{\centering\arraybackslash}p{0.475cm}>{\centering\arraybackslash}p{0.475cm}>{\centering\arraybackslash}p{0.475cm}>{\centering\arraybackslash}p{0.475cm}>{\centering\arraybackslash}p{0.475cm}>{\centering\arraybackslash}p{0.475cm}>{\centering\arraybackslash}p{0.475cm}}
$n=$ &  1 &  2 &  3 &  4 &  5 &  6 &  7 &  8 &  9 & 10 & 11 & 12 & 13 & 14 & 15 & 16 \\ \hline
Linear      &0&  1&  3&  6& 10& 15& 21& 28& 36& 45& 55& 66& 78& 91&105&120 \\
Shellsort   &0&  1&  3&  6& 10& 14& 18& 23& 29& 35& 41& 48& 56& 64& 71& 78 \\
\hline
Improvement &0&  0&  0&  0&  0&  1&  3&  5&  7& 10& 14& 18& 22& 27& 34& 42 \\
\end{tabular}
\caption{Minimum Numbers of Comparisons to Guarantee a Complete Sort by Different Methods}
\label{tab:Shell_Linear}
\end{table}

\vspace*{\fill}\begin{table}[H]
\begin{minipage}[t]{.5\linewidth}\centering
\begin{tabular}{p{1.5cm}p{3.5cm}p{1.5cm}}
$n=1$ &  & \\ \hline
terminated. &  & \\
    &  & \\
$n=2$ &  & \\ \hline
$i=1$ & $s_i= 1 $  & $n_i=1$ \\
terminated. &  & \\
    &  & \\
$n=3$ &  & \\ \hline
$i=1$ & $s_i= 1 $  & $n_i=3$ \\
terminated. &  & \\
    &  & \\
$n=4$ &  & \\ \hline
$i=1$ & $s_i= 1 $  & $n_i=6$ \\
terminated. &  & \\
    &  & \\
$n=5$ &  & \\ \hline
$i=1$ & $s_i= 1 $  & $n_i=10$ \\
terminated. &  & \\
    &  & \\
$n=6$ &  & \\ \hline
$i=1$ & $s_i= 1 $  & $n_i=15$ \\
$i=5$ & $s_i= 1,\, 4 $  & $n_i=14$ \\
terminated. &  & \\
    &  & \\
$n=7$ &  & \\ \hline
$i=1$ & $s_i= 1$  & $n_i=21$ \\
$i=3$ & $s_i= 1,\, 3$  & $n_i=20$ \\
$i=5$ & $s_i= 1,\, 4$  & $n_i=19$ \\
$i=21$ & $s_i= 1,\, 4,\, 6$  & $n_i=18$ \\
terminated. &  & \\
    &  & \\
$n=8$ &  & \\ \hline
$i=1$ & $s_i= 1 $  & $n_i=28$ \\
$i=3$ & $s_i= 1,\, 3 $  & $n_i=25$ \\
$i=7$ & $s_i= 1,\, 3,\, 4 $  & $n_i=24$ \\
$i=41$ & $s_i= 1,\, 5,\, 7 $  & $n_i=23$ \\
terminated. &  &   
\end{tabular}
\end{minipage}\hfill
\begin{minipage}[t]{.5\linewidth}\centering
\centering
\begin{tabular}{p{1.5cm}p{3.5cm}p{1.5cm}}
$n=9$ &  & \\ \hline
$i=1$ & $s_i= 1 $  & $n_i=36$ \\
$i=2$ & $s_i= 1,\, 2$  & $n_i=34$ \\
$i=3$ & $s_i= 1,\, 3$  & $n_i=33$ \\
$i=4$ & $s_i= 1,\, 2,\, 3$  & $n_i=32$ \\
$i=7$ & $s_i= 1,\, 3,\, 4 $  & $n_i=29$ \\
terminated. &  & \\
    &  & \\
$n=10$ &  & \\ \hline
$i=1$ & $s_i= 1 $  & $n_i=45$ \\
$i=2$ & $s_i= 1,\, 2$  & $n_i=43$ \\
$i=3$ & $s_i= 1,\, 3$  & $n_i=39$ \\
$i=5$ & $s_i= 1,\, 4$  & $n_i=37$ \\
$i=7$ & $s_i= 1,\, 3,\, 4 $  & $n_i=36$ \\
$i=145$ & $s_i= 1,\, 6,\, 9 $  & $n_i=35$ \\
terminated. &  & \\
    &  & \\
$n=11$ &  & \\ \hline
$i=1$ & $s_i= 1 $  & $n_i=55$ \\
$i=2$ & $s_i= 1,\, 2 $  & $n_i=50$ \\
$i=3$ & $s_i= 1,\, 3 $  & $n_i=46$ \\
$i=5$ & $s_i= 1,\, 4 $  & $n_i=45$ \\
$i=7$ & $s_i= 1,\,3,\, 4 $  & $n_i=43$ \\
$i=13$ & $s_i= 1,\, 4,\, 5 $  & $n_i=41$ \\
terminated. &  & \\
    &  & \\
$n=12$ &  & \\ \hline
$i=1$ & $s_i= 1 $  & $n_i=66$ \\
$i=2$ & $s_i= 1,\, 2 $  & $n_i=61$ \\
$i=3$ & $s_i= 1,\, 3 $  & $n_i=57$ \\
$i=4$ & $s_i= 1,\, 2,\, 3 $  & $n_i=53$ \\
$i=7$ & $s_i= 1,\, 3,\, 4 $  & $n_i=49$ \\
$i=547$ & $s_i= 1,\, 3,\, 7,\, 11 $  & $n_i=48$ \\
terminated. &  & \\
 &  & \\
 &  & \\
 &  &
 \end{tabular}
\end{minipage}
\caption{Complete Search $n=1,2,3,4,5,6,7,8,9,10,11,12$}
\label{table:12}
\end{table}\vspace*{\fill}

\vspace*{\fill}\begin{table}[H]
\begin{minipage}[t]{.5\linewidth}\centering
\begin{tabular}{p{1.5cm}p{3.5cm}p{1.5cm}}
$n=13$ &  & \\ \hline
$i=1$ & $s_i= 1$  & $n_i=78$ \\
$i=2$ & $s_i= 1,\, 2 $  & $n_i=69$ \\
$i=3$ & $s_i= 1,\, 3 $  & $n_i=64$ \\
$i=4$ & $s_i= 1,\, 2,\, 3 $  & $n_i=61$ \\
$i=7$ & $s_i= 1,\, 3,\, 4 $  & $n_i=56$ \\
terminated. &  & \\
    &  & \\
$n=14$ &  & \\ \hline
$i=1$ & $s_i= 1 $  &       $n_i=91$ \\
$i=2$ & $s_i= 1,\, 2 $     & $n_i=82$ \\
$i=3$ & $s_i= 1,\, 3 $     & $n_i=73$ \\
$i=4$ & $s_i= 1,\, 2,\, 3 $  & $n_i=70$ \\
$i=7$ & $s_i= 1,\, 3,\, 4 $  & $n_i=64$ \\
terminated. &  & \\
    &  & \\
$n=15$ &  & \\ \hline
$i=1$ & $s_i= 1$  & $n_i=105$ \\
$i=2$ & $s_i= 1,\, 2$  & $n_i=91$ \\
$i=3$ & $s_i= 1,\, 3$  & $n_i=87$ \\
$i=4$ & $s_i= 1,\, 2, \,3$  & $n_i=80$ \\
$i=7$ & $s_i= 1,\, 3, \,4$  & $n_i=72$ \\
$i=35$ & $s_i= 1,\, 3, \,7$  & $n_i=71$ \\
terminated. &  & \\
    &  & \\
$n=16$ &  & \\ \hline
$i=1$ & $s_i= 1$  & $n_i=120$ \\
$i=2$ & $s_i= 1,\, 2$  & $n_i=106$ \\
$i=3$ & $s_i= 1,\, 3$  & $n_i=95$ \\
$i=4$ & $s_i= 1,\, 2,\, 3$  & $n_i=89$ \\
$i=7$ & $s_i= 1,\, 3,\, 4$  & $n_i=79$ \\
$i=165$ & $s_i= 1,\, 4,\, 7,\, 9$  & $n_i=78$ \\
terminated. &  &  
\end{tabular}\caption{Complete Search $n=13,14,15,16$}
\label{table:16}\end{minipage}\hfill
\begin{minipage}[t]{.5\linewidth}\centering
\centering
\begin{tabular}{p{1.5cm}p{3.5cm}p{1.5cm}}
$n=17$ &  & \\ \hline
$i=1$ & $s_i= 1 $  & $n_i=136$ \\
$i=2$ & $s_i= 1,\, 2 $  & $n_i=116$ \\
$i=3$ & $s_i= 1,\, 3 $  & $n_i=106$ \\
$i=4$ & $s_i= 1,\, 2,\, 3 $  & $n_i=101$ \\
$i=7$ & $s_i= 1,\, 3,\, 4 $  & $n_i=87$ \\
$i=45$ & $s_i= 1,\, 4,\, 5,\, 7 $  & $n_i\geq86$ \\
    &  & \\
$n=18$ &  & \\ \hline
$i=1$ & $s_i= 1 $  &       $n_i=153$ \\
$i=2$ & $s_i= 1,\, 2 $     & $n_i=133$ \\
$i=3$ & $s_i= 1,\, 3 $     & $n_i=123$ \\
$i=4$ & $s_i= 1,\, 2,\, 3 $  & $n_i=109$ \\
$i=7$ & $s_i= 1,\, 3,\, 4 $  & $n_i=100$ \\
$i=10$ & $s_i= 1,\, 2,\, 5 $  & $n_i= 98$ \\
$i=34$ & $s_i= 1,\, 2,\, 7 $  & $n_i\geq 97$ \\
    &  & \\
$n=19$ &  & \\ \hline
$i=1$ & $s_i= 1 $  &       $n_i=171$ \\
$i=2$ & $s_i= 1,\, 2 $     & $n_i=144$ \\
$i=3$ & $s_i= 1,\, 3 $     & $n_i=132$ \\
$i=4$ & $s_i= 1,\, 2,\, 3 $  & $n_i=119$ \\
$i=7$ & $s_i= 1,\, 3,\, 4 $  & $n_i=109$ \\
$i=11$ & $s_i= 1,\, 3,\, 5 $  & $n_i=105$ \\
$i=45$ & $s_i= 1,\, 4,\, 5,\, 7$  & $n_i\geq 104$ \\
    &  & \\
$n=20$ &  & \\ \hline
$i=1$ & $s_i= 1 $  &       $n_i=190$ \\
$i=2$ & $s_i= 1,\, 2 $     & $n_i=163$ \\
$i=3$ & $s_i= 1,\, 3 $     & $n_i=145$ \\
$i=4$ & $s_i= 1,\, 2,\, 3 $  & $n_i=137$ \\
$i=7$ & $s_i= 1,\, 3,\, 4 $  & $n_i=117$ \\
$i=11$ & $s_i= 1,\, 3,\, 5 $  & $n_i\geq 113$  
\end{tabular}
\caption{Incomplete Search $n=17,18,19,20$}
\label{table:20}
\end{minipage}
\end{table}\vspace*{\fill}

\clearpage\vspace*{\fill}\begin{table}[H]
\begin{minipage}[t]{.5\linewidth}\centering
\begin{tabular}{p{1.5cm}p{3.5cm}p{1.5cm}}
$n=21$ &  & \\ \hline
$i=1$ & $s_i= 1$                & $n_i=210$ \\
$i=2$ & $s_i= 1,\, 2 $          & $n_i=175$ \\
$i=3$ & $s_i= 1,\, 3 $          & $n_i=165$ \\
$i=4$ & $s_i= 1,\, 2,\, 3 $     & $n_i=146$ \\
$i=7$ & $s_i= 1,\, 3,\, 4 $     & $n_i=126$ \\
$i=11$ & $s_i= 1,\, 3,\, 5 $     & $n_i\geq 122$ \\
    &  & \\
$n=22$ &  & \\ \hline
$i=1$ & $s_i= 1 $  &              $n_i=231$ \\
$i=2$ & $s_i= 1,\, 2 $          & $n_i=196$ \\
$i=3$ & $s_i= 1,\, 3 $          & $n_i=175$ \\
$i=4$ & $s_i= 1,\, 2,\, 3 $     & $n_i=157$ \\
$i=7$ & $s_i= 1,\, 3,\, 4 $     & $n_i\geq 140$ \\
    &  & \\
$n=23$ &  & \\ \hline
$i=1$ & $s_i= 1 $  &              $n_i=253$ \\
$i=2$ & $s_i= 1,\, 2 $          & $n_i=209$ \\
$i=3$ & $s_i= 1,\, 3 $          & $n_i=190$ \\
$i=4$ & $s_i= 1,\, 2,\, 3 $     & $n_i=173$ \\
$i=7$ & $s_i= 1,\, 3,\, 4 $     & $n_i\geq 145$ \\
    &  & \\
$n=24$ &  & \\ \hline
$i=1$ & $s_i= 1 $  &              $n_i=276$ \\
$i=2$ & $s_i= 1,\, 2 $          & $n_i=232$ \\
$i=3$ & $s_i= 1,\, 3 $          & $n_i=213$ \\
$i=4$ & $s_i= 1,\, 2,\, 3 $     & $n_i=183$ \\
$i=7$ & $s_i= 1,\, 3,\, 4 $     & $n_i\geq 153$ \\
    &  & \\
$n=25$ &  & \\ \hline
$i=1$ & $s_i= 1 $  &              $n_i=300$ \\
$i=2$ & $s_i= 1,\, 2 $          & $n_i=246$ \\
$i=3$ & $s_i= 1,\, 3 $          & $n_i=224$ \\
$i=4$ & $s_i= 1,\, 2,\, 3 $     & $n_i=195$ \\
$i=7$ & $s_i= 1,\, 3,\, 4 $     & $n_i\geq 161$ 
\end{tabular}\end{minipage}\hfill
\begin{minipage}[t]{.5\linewidth}\centering
\centering
\begin{tabular}{p{1.5cm}p{3.5cm}p{1.5cm}}
$n=26$ &  & \\ \hline
$i=1$ & $s_i= 1 $  &              $n_i=325$ \\
$i=2$ & $s_i= 1,\, 2 $          & $n_i=271$ \\
$i=3$ & $s_i= 1,\, 3 $          & $n_i=241$ \\
$i=4$ & $s_i= 1,\, 2,\, 3 $     & $n_i=219$ \\
$i=7$ & $s_i= 1,\, 3,\, 4 $     & $n_i\geq 174$ \\
    &  & \\
$n=27$ &  & \\ \hline
$i=1$ & $s_i= 1 $  &              $n_i=351$ \\
$i=2$ & $s_i= 1,\, 2 $          & $n_i=286$ \\
$i=3$ & $s_i= 1,\, 3 $          & $n_i=267$ \\
$i=4$ & $s_i= 1,\, 2,\, 3 $     & $n_i=230$ \\
$i=7$ & $s_i= 1,\, 3,\, 4 $     & $n_i\geq 188$ \\
    &  & \\
$n=28$ &  & \\ \hline
$i=1$ & $s_i= 1 $  &              $n_i=378$ \\
$i=2$ & $s_i= 1,\, 2 $          & $n_i=313$ \\
$i=3$ & $s_i= 1,\, 3 $          & $n_i=279$ \\
$i=4$ & $s_i= 1,\, 2,\, 3 $     & $n_i=243$ \\
$i=7$ & $s_i= 1,\, 3,\, 4 $     & $n_i\geq 193$ \\
    &  & \\
$n=29$ &  & \\ \hline
$i=1$ & $s_i= 1 $  &              $n_i=406$ \\
$i=2$ & $s_i= 1,\, 2 $          & $n_i=329$ \\
$i=3$ & $s_i= 1,\, 3 $          & $n_i=298$ \\
$i=4$ & $s_i= 1,\, 2,\, 3 $     & $n_i=263$ \\
$i=7$ & $s_i= 1,\, 3,\, 4 $     & $n_i\geq 212$ \\
    &  & \\
$n=30$ &  & \\ \hline
$i=1$ & $s_i= 1 $  &              $n_i=435$ \\
$i=2$ & $s_i= 1,\, 2 $          & $n_i=358$ \\
$i=3$ & $s_i= 1,\, 3 $          & $n_i=327$ \\
$i=4$ & $s_i= 1,\, 2,\, 3 $     & $n_i=275$ \\
$i=7$ & $s_i= 1,\, 3,\, 4 $     & $n_i\geq 224$ \\
    &  & 
\end{tabular}
\end{minipage}
\caption{Incomplete Search $n=21,22,23,24,25,26,27,28,29,30$}
\label{table:30}
\end{table}\vspace*{\fill}

\clearpage\section{Summary}
Define for any $i\in\mathbb{N}_1$, $n_i:\mathbb{N}_1\to\mathbb{N}_0$ by for any $n\in\mathbb{N}_1$, $n_i(n)$ to be the maximum number of comparisons required to sort $n$ pairwise district elements by the gap sequence $s_i$ in Shellsort. Define for any $k\in\mathbb{N}_1$, $\chi_k:\mathbb{N}_1\to\{0,1\}$ by for any $n\in\mathbb{N}_1$, $\chi_k(n)=1$ if and only if $n$ is an integer multiple of $k$.

From the search history, for any $n\in\mathbb{N}_3$, the explicit formula of $n_2(n)$ is given by \begin{align}
n_2(n)
&=\frac{n(n-1)}{2}-\frac12\left\lceil\frac{n}{2}\right\rceil^2+\frac52\left\lceil\frac{n}{2}\right\rceil+2\\
&=-\frac{7}{8}+\frac{n}{2}+\frac{3n^2}{8}-\frac{9}{8}\chi_2(n)+\frac{n}{4}\chi_2(n)\end{align}

$n_2(n)<n(n-1)/2$, the maximum number of comparisons required to sort $n$ pairwise district elements by linear insertion sort, if and only if $n\geq9$. Similarly, for any $n\in\mathbb{N}_4$, \begin{align}
n_3(n)
&=\frac{n(n-1)}{2}-\frac32\left\lceil\frac{n}{3}\right\rceil^2+\frac92\left\lceil\frac{n}{3}\right\rceil-3+\frac{n-1}{3}\left(\left\lceil\frac{n}{3}\right\rceil-\left\lceil\frac{n-1}{3}\right\rceil\right)\\
&=-1+\frac{2n}{3}+\frac{n^2}{3}-\left(2-\frac{n}{3}\right)\chi_3(n)-\frac{2}{3}\chi_3(n+1)\end{align}

$n_3(n)<n(n-1)/2$ if and only if $n\geq7$, and $n_3(n)<n_2(n)$ if and only if $n\geq9$. For any $n\in\mathbb{N}_4$, \begin{multline}
n_4(n)
= -\frac{35}{12} + \frac{5n}{3} + \frac{n^2}{4}-\left(\frac{25}{12}-\frac{n}{6}\right)\chi_6(n)-\left(\frac{7}{3}-\frac{n}{3}\right)\chi_6(n+1)\\-\left(\frac{17}{12}-\frac{n}{6}\right)\chi_6(n+2)-\left(\frac{10}{3}-\frac{n}{3}\right)\chi_6(n+3)-\left(\frac{41}{12}-\frac{n}{2}\right)\chi_6(n+4)\end{multline}

Although many further similar inequalities are suggested from the search history, two of the most notable that describing for long term behaviour are descending sequence that for any $n\in\mathbb{N}_9$,
\begin{align}\label{n7}
n_7(n)<n_4(n)<n_3(n)<n_2(n)<n_1(n)\end{align}and that for any $n\in\mathbb{N}_{19}$,
\begin{align}\label{n11}
n_{11}(n)<n_7(n)<n_4(n)<n_3(n)<n_2(n)<n_1(n)\end{align}

\vfill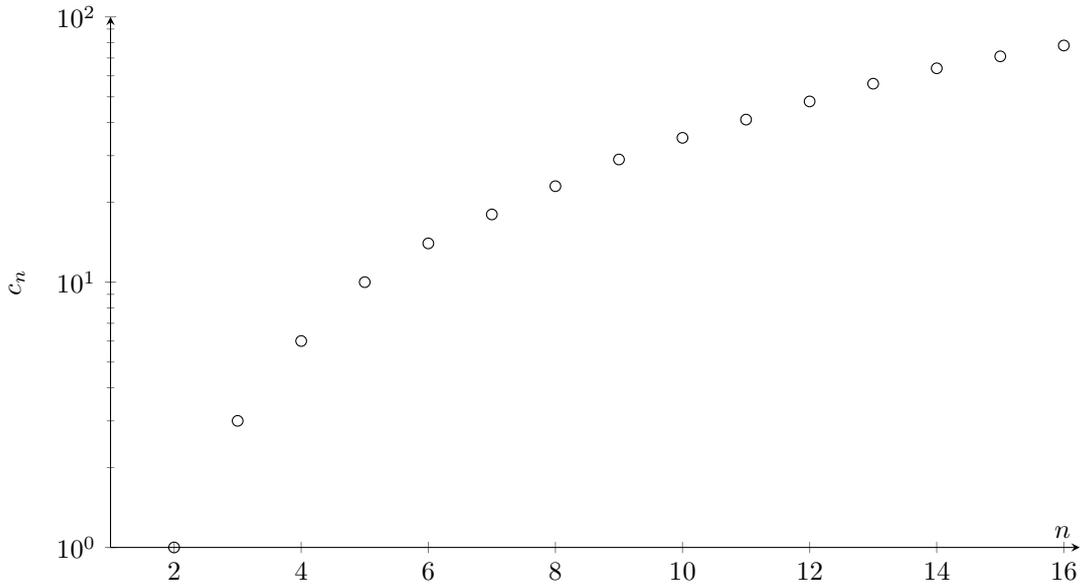
\begin{figure}[h!]\begin{tikzpicture}\centering
\begin{axis}[axis lines = left,
 width=0.9\textwidth,
 height=0.35\textheight,
 axis x line=middle,
 hide obscured x ticks=false,
 xlabel = $n$,
 ylabel = $c_n$,
 xmode=normal,ymode=log,
 xmin=1, xmax=16.25,
 ymin=0, ymax=100,
]
\addplot+[only marks,mark=o,draw=black]coordinates {
( 2, 1) ( 3, 3) ( 4, 6) ( 5,10)
( 6,14) ( 7,18) ( 8,23) ( 9,29)
(10,35) (11,41) (12,48) (13,56)
(14,64) (15,71) (16,78) };\end{axis}
\end{tikzpicture}\caption{Minimax Number of Comparisons $c_n$ to Sort $n$ Pairwise Distinct Elements in Shellsort}\label{2}\end{figure}

\clearpage\section{Conclusion}
Similarly, for any $n\in\mathbb{N}_1$ and $s\in S_n$, the set of all Bad $(s,2)$-sorted permutations is defined to be the set of all Bad $(s,1)$-sorted permutations that are Bad $(s(2),1)$-sorted permutations after it completes the first pass of $s$, where $s(2)$ denotes the second largest element in $s$. In mathematical notation, the set of all Bad $(s,2)$-sorted permutations is defined as\begin{align}
P_{n,{(s,2)}}=\left\{p\in P_{n,{(s,1)}}:p_{s(1)}\in P_{n,{(s(2),1)}}\right\}\end{align}
, where $p_{s(1)}$ denotes the permutation of $p$ after it completes the first pass of $s$. However, most often it is the empty set, for example when $n=16$ and $s=\{1,4\}$, $P_{n,{(s,2)}}=\emptyset$.

Further simplification for the minimax search seems to be difficult in this aspect, instead, reordering the index $i$ for the gap sequence in the minimising step, as well as reordering the elements in $P_{n,{(s,1)}}$ for $s$ to sort, according to the value of $s(2)$, in the maximising step seems to be more practical. 

Best average comparison gap sequences can be defined and obtained, in a similar way to optimal gap sequences, in the sense of minimised average number of comparisons. Simplified algorithm to establish those sequences seems to be hard to implement, but only being searched by brute-force method. 

The optimal gap sequence can also be defined to have further more properties on top of having the minimised maximum number of comparisons, such as in addition to having the minimised maximum number of swaps, having the minimum average number of comparisons or containing the least number of increments, among all optimal gap sequences.  

\vfill
\begin{figure}[h]
\begin{subfigure}[b]{0.49\linewidth}
\centering
\begin{tikzpicture}
\begin{axis} [width=1.15\textwidth,height=0.34\textheight, title={Linear Insertion Sort ($s_1=\{1\}$)}, ybar,bar width=16pt,clip=false,axis y line=none,axis x line*=bottom,nodes near coords,ymin=0,ymax=210,xtick = {5,...,15}]]
\addplot [fill=none,xtick=data]coordinates { 
 ( 5,  2)
 ( 6, 10) 
 ( 7, 26) 
 ( 8, 52)
 ( 9, 82) 
 (10,110) 
 (11,126) 
 (12,120)
 (13, 96) 
 (14, 64) 
 (15, 32)
};
\end{axis} 
\end{tikzpicture}
\end{subfigure}
\begin{subfigure}[b]{0.49\linewidth}
\centering
\begin{tikzpicture}
\begin{axis} [width=1.15\textwidth,height=0.34\textheight,title={Optimal Gap Sequence ($s^{**}_6=s_5=\{1,\,4\}$) }, ybar,bar width=16pt,clip=false,axis y line=none,axis x line*=bottom,nodes near coords,ymin=0,ymax=210,xtick = {5,...,15}]
\addplot [fill=none,xtick=data]coordinates { 
 ( 5,  0)
 ( 6,  0) 
 ( 7,  8) 
 ( 8, 40)
 ( 9,104) 
 (10,180) 
 (11,192) 
 (12,128)
 (13, 56) 
 (14, 12) 
 (15,  0)
};
\end{axis} 
\end{tikzpicture}\end{subfigure}\caption{Distributions of the Numbers of Comparisons to Sort $n=6$ Pairwise Distinct Elements}\end{figure}
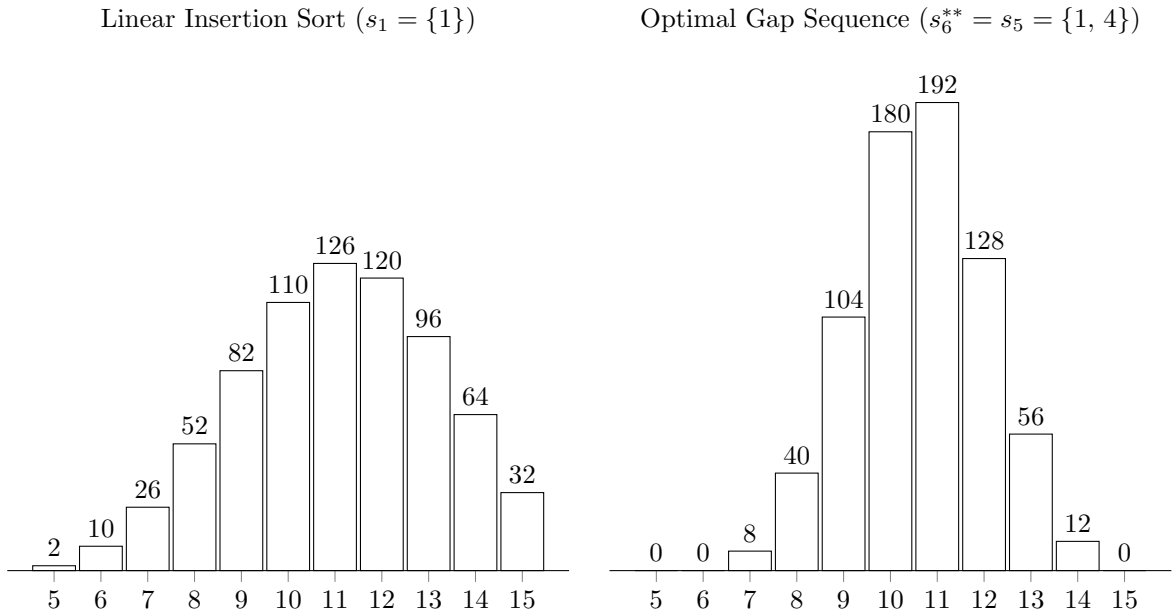

\end{document}